\documentclass{amsart}
\usepackage{amsmath,amssymb}
\usepackage[latin1]{inputenc}
\newtheorem{theorem}{Theorem}[section]

\newtheorem{corollary}[theorem]{Corollary}

\newtheorem{lemma}[theorem]{Lemma}

\theoremstyle{definition}

\theoremstyle{remark}
\newtheorem{remark}[theorem]{Remark}

\numberwithin{equation}{section}

\newcommand{\al}{\alpha}
\newcommand{\be}{\beta}
\newcommand{\de}{\delta}
\newcommand{\ep}{\epsilon}

\newcommand{\ga}{\gamma}

\newcommand{\la}{\lambda}

\newcommand{\si}{\sigma}

\newcommand{\vp}{\varphi}

\newcommand{\De}{\Delta}
\newcommand{\Ga}{\Gamma}
\newcommand{\La}{\Lambda}
\newcommand{\Si}{\Sigma}

\def\NN{\mathbb{N}}
\def\RR{\mathbb{R}}

\newcommand{\cG}{{\mathcal G}}

\newcommand{\cK}{{\mathcal K}}

\newcommand{\cS}{{\mathcal S}}

\newcommand{\pd}{\partial}
\newcommand\minus\backslash

\newcommand\lan\langle
\newcommand\ran\rangle

\newcommand{\tr}{\operatorname{tr}}
\newcommand{\im}{\operatorname{im}}

\newcommand{\supp}{\operatorname{supp}}

\renewcommand\leq\leqslant
\renewcommand\geq\geqslant
\newlength{\intwidth}

\DeclareMathOperator\ev{ev}

 \DeclareMathOperator\ind{ind}
\DeclareMathOperator\codim{codim}

\begin{document}

\title[Nondegeneracy of the eigenvalues of the Hodge Laplacian]{Nondegeneracy of the eigenvalues of the Hodge Laplacian for generic metrics on
  3-manifolds}

\author{Alberto Enciso}
\address{Departement Mathematik, ETH Zürich, 8092 Zürich, Switzerland}
\email{alberto.enciso@math.ethz.ch}

\author{Daniel Peralta-Salas}
\address{Instituto de Ciencias Matemáticas, CSIC-UAM-UC3M-UCM, C/ Serrano 123, 28006 Madrid, Spain. Current address: Departamento de Matemáticas, Universidad Carlos III, 28911 Leganés, Spain}
\email{dperalta@icmat.es}

%
%
\begin{abstract}
  In this paper we analyze the eigenvalues and eigenfunctions of the Hodge Laplacian for generic metrics on a closed 3-manifold $M$. In particular, we show that the nonzero eigenvalues are simple and the zero set of the eigenforms of degree $1$ or $2$ consists of isolated points for a residual set of $C^r$ metrics on $M$, for any integer $r\geq2$. The proof of this result hinges on a detailed study of the Beltrami (or rotational) operator on co-exact $1$-forms. 
\end{abstract}
\maketitle

\section{Introduction and statement of results}
\label{S:intro}

The analysis of the eigenvalues of the Laplacian on a
compact manifold is a multifaceted subject of central interest in
Riemannian geometry. From a qualitative point of view, one of the most beautiful results in this regard is a celebrated theorem of
Uhlenbeck~\cite{Uh76} which ensures that, given a closed
manifold $M$, the set of $C^r$ metrics on $M$ whose Laplacian has simple
spectrum is residual, for any $2\leq r<\infty$.

Shortly after the publication of~\cite{Uh76}, it
was observed by Millman~\cite{Mi80} that Uhlenbeck's theorem cannot
hold true for the Hodge Laplacian without further hypotheses. Indeed,
if $M$ has dimension $2n$, Millman showed that any eigenvalue of the
Laplacian on $n$-forms must have even multiplicity by the
McKean--Singer télescopage theorem. The question of whether an appropriate analog of
Uhlenbeck's theorem is valid for the Hodge Laplacian was left wide
open, and is actually included in the problem section of~\cite[Problem 8.24]{CPR01}.

The study of the zero set of the eigenfunctions of the Laplacian has also attracted considerable attention. It has long been known~\cite{HS89} that the zero set of a scalar eigenfunction of the Laplacian  in an $n$-manifold is a countably $(n-1)$-rectifiable set with finite $(n-1)$-Hausdorff measure, and that it is in fact a codimension~$1$ submanifold of class  $C^{r,\al}$ for a residual set of $C^r$ metrics~\cite{Uh76}. On the contrary, our understanding of the zero set of the eigenfunctions of the Laplacian on $p$-forms is rather limited, for $p\neq 0,n$. Since the zero set of an exact $1$-form $df$ which is an eigenfunction of the Laplacian is simply the critical set of the corresponding scalar eigenfunction, we do know that it consists of isolated points for a residual set of $C^r$ metrics~\cite{Uh76}. In addition to this,  a theorem of Bär~\cite{Ba99} ensures that the zero set of a nontrivial harmonic form has finite $(n-2)$-measure, but otherwise the zero set of an eigenform of positive eigenvalue can have positive $(n-1)$-measure.

Our purpose in this paper is to show that, given a closed three-dimensional manifold $M$, there exists a residual set of $C^r$ metrics such that all the nonzero eigenvalues of the Hodge Laplacian on $p$-forms have multiplicity $1$, for all $0\leq p\leq3$. Moreover, we also prove that the zero set of the corresponding eigenfunctions consists of isolated points for $p=1,2$. The fact that the zero set of the eigenforms should have codimension higher than $1$ in some generic sense had been conjectured by Yau in~\cite[Problem 38]{Ya93}; in dimension $3$, we show that this codimension is in fact maximal. It would be of great interest to ascertain whether the cardinality of the zero set is generically increasing with the eigenvalue. For exact $1$-forms, this is equivalent to the question posed in~\cite[Problem 76]{Ya82}, which one can hope to answer in the affirmative at most for generic metrics in view of Jakobson and Nadirashvili's example of metrics admitting an infinite sequence of eigenfunctions with a uniformly bounded number of critical points~\cite{JN99}. 

Our main result can be stated as follows:
  
\begin{theorem}\label{T.main}
  Given a closed $3$-manifold $M$ and an integer $r\geq2$, there exists a residual subset $\Ga$ of the space of $C^r$ metrics on $M$ such that, for all $g\in \Ga$, the nonzero eigenvalues of the Hodge Laplacian $\De_g$ on $p$-forms have multiplicity $1$ for all $0\leq p\leq3$ and the zero set of each eigenfunction consists of isolated points for $p=1,2$. More precisely,
  \begin{enumerate}
  \item The nonzero eigenvalues of $\De_g$ on $1$-forms are simple and the corresponding eigenfunctions have all hyperbolic zeros.
  \item The eigenvalues of $\De_g$ on functions are simple and the corresponding eigenfunctions are Morse.
  \item The nonzero eigenvalues of $\De_g$ on co-closed $1$-forms and on functions are disjoint.
  \end{enumerate}
\end{theorem}

Some comments on the strategy of the proof are in order. A first observation is that the cases of $0$-forms and of exact $1$-forms can be readily dealt with using Uhlenbeck's theorem for the scalar Laplacian. By Hodge duality, the crux of the matter is then to prove a generic nondegeneracy result for the Laplacian on co-exact $1$-forms via infinite-dimensional transversality theory. A formula for the variation of simple eigenvalues is ultimately used to separate the spectrum of the Laplacian on exact and co-exact $1$-forms.

The topological obstruction to simplicity of eigenvalues observed by Millman shows that fundamental subtleties appear in the extension of Uhlenbeck's theorem to differential forms. Therefore, in our study of the Laplacian on co-exact $1$-forms we follow an indirect approach which hinges on the analysis of the spectrum of the Beltrami (or rotational) operator, thus making essential use of the fact that the manifold is three-dimensional. In spite of the fact that the Beltrami operator is simpler than the Hodge Laplacian, various technical difficulties make the extension of Uhlenbeck's proof to differential forms rather nontrivial. These are mainly due to the facts that the space of co-exact $1$-forms depends on the choice of the metric, that the symbol of the Beltrami operator is not elliptic and that the PDE is now vector-valued. A thorough discussion of these points is included in Section~\ref{S.Beltrami}; in a way, they reflect many of the difficulties which make the spectral theoretic analysis of the Hodge Laplacian rather different than its scalar counterpart, and which account for the substantial recent interest attracted by the qualitative properties of the differential form spectrum~\cite{AC95,CT97,Lo02,Ma08,Ja09}.

As a side remark, let us notice that the analysis of the Hodge Laplacian on a surface is much less involved than in three dimensions. In fact, an easy application of Uhlenbeck's theorem shows that, for a residual set of metrics, the nonzero eigenvalues of the Laplacian on $p$-forms are simple for $p=0,2$ and have multiplicity $2$ for $p=1$.

The article is organized as follows. In Section~\ref{S.lemmas}, we prove several technical lemmas that will be required in later sections. In Section~\ref{S.Beltrami}, the generic properties of the eigenvalues and eigenfunctions of the Beltrami operator are addressed using the machinery of infinite-dimensional transversality theory. A brief discussion of the new technical complications that arise in the extension of Uhlenbeck's theorem to the Beltrami operator is included. In Section~\ref{S.Hodge}, we provide the proof of Theorem~\ref{T.main}, which makes use of all the previous results derived in this article.

\section{Some preliminary lemmas}
\label{S.lemmas}

Throughout this article, $M$ will stand for an oriented compact
3-manifold of class $C^\infty$ without boundary. We shall consider the space $\cG^r(M)$ of metrics on $M$ of class $C^r$, for some integer $r\geq2$. $\cG^r(M)$ is then a $C^r$ Banach manifold, whose tangent space at any point can be identified with the space $\cS^{r}(M)$ of symmetric tensor fields of class $C^{r}$ and type $(0,2)$.

The {\em Beltrami operator} (or rotational) associated to the metric $g$ is $*_gd$, acting on $1$-forms. By $d$, $\de_g$ and $*_g$ we shall respectively denote the
differential, codifferential and star operator computed with respect to the metric $g$. The spaces of $C^k$, $L^2$ and $H^k$ $p$-forms in $M$ will be respectively denoted by
$C^k(M,\La^p)$, $L^2(M,\La^p)$ and $H^k(M,\La^p)$. The $L^2$ inner product associated to $g$ is
\begin{equation}\label{lanran}
\lan u,v\ran_g:=\int u\wedge *_g v\,,
\end{equation}
and the associated norm is denoted by $\|\cdot\|_g$. As a topological vector space, the Sobolev space $H^k(M,\La^p)$ is obviously independent of the $C^r$ Riemannian metric on $M$ we use to define the $H^k$ scalar product, for all $0\leq k\leq r$.

Let us consider the closed subspace
\[
\cK:=\big\{ u\in L^2(M,\La^1):du=0\big\}
\]
of $L^2(M,\La^1)$ and define the fiber bundle
\[
E:=\big\{(g,u):g\in\cG^r(M),\; u\in E_g\big\}
\]
over $\cG^r(M)$ whose fiber at $g$ is
\begin{equation}\label{Eg}
E_g:=\big\{u\in H^1(M,\La^1): u\perp_g\cK\,,\; \|u\|_g=1\big\}\,.
\end{equation}
The symbol $\perp_g$ denotes orthogonality with respect to the inner product~\eqref{lanran}, and we  will use the notation $S^{\perp_g}$ for the orthogonal complement of a closed subspace $S\subset L^2(M,\La^1)$ with respect to~\eqref{lanran}. It is well known that $*_gd$ defines an unbounded self-adjoint operator on the Hilbert space $(\cK^{\perp_g},\lan\cdot,\cdot\ran_g)$ with domain $ H^1(M,\La^1)\cap \cK^{\perp_g}$, which we shall also call Beltrami operator.

As customary, we respectively denote by $g_{ij}$ and $g^{ij}$ the components of the metric tensor  in local coordinates and its inverse matrix. The Riemannian volume $3$-form is denoted by $\mu_g$, and with a slight abuse of notation we will call $d\mu_g$ the corresponding volume measure. The metric will be used to raise the indices of covariant tensor fields; e.g., for any tensor fields of components $w_i$ and $T_{ij}$ we define
\begin{equation}\label{raise}
w^i:=g^{ij}w_j\,,\qquad T^{ij}:=g^{ik}g^{jl}T_{kl}\,.
\end{equation}
(The metric used to raise the indices will be clear from the context). The trace of a tensor $h\in\cS^r(M)$ is $\tr_gh:=g^{ij}h_{ij}$.

Let us begin by evaluating the variation of the Beltrami operator with respect to the metric. We  use the customary notation $(D\psi)_x:T_xX\to T_{\psi(x)}Y$ for the derivative at $x$ of a differentiable map $\psi:X\to Y$, so that by $D(*d)_g(h)$ we will denote the variation of the Beltrami operator at the metric $g\in\cG^r(M)$ in the direction determined by the tensor $h\in\cS^r(M)$.

\begin{lemma}\label{L.D}
  Let $u\in E_g$ be an eigenfunction of $*_gd$ with eigenvalue $\la$. Then, for any $h\in\cS^r(M)$,
  \[
\big(D(*d)_g(h)u\big)_i=\la\, h_{ij}u^j-\frac\la2\,(\tr_g h)\, u_i\,.
  \]
\end{lemma}
\begin{proof}
Notice that
\[
\big(*_gdu\big)_k=\frac12(\mu_g)_{ijk}g^{il}g^{jm}\bigg(\frac{\pd u_m}{\pd x^l}-\frac{\pd u_l}{\pd x^m}\bigg)\,,
\]
where $(\mu_g)_{ijk}=|g|^{1/2}\ep_{ijk}$, $|g|:=\det(g_{ij})$ and  $\ep_{ijk}$ is the Levi-Civita permutation symbol. Taking derivatives in the above formula and using that $D(|g|^s)(h)=s|g|^s\tr_gh$ and $D(g^{ij})(h)=-h^{ij}$, where $h^{ij}$ is defined as in~\eqref{raise} and $s>0$, one immediately derives that
\begin{equation}\label{DG*du}
\big(D(*d)_g(h)u\big)_k=  \frac12(\mu_g)_{ijk} \bigg(\frac12\tr_g h\,g^{il}g^{jm} - h^{il}g^{jm}-g^{il}h^{jm}\bigg)\bigg(\frac{\pd u_m}{\pd x^l}-\frac{\pd u_l}{\pd x^m}\bigg)\,.
\end{equation}
The equation $*_gdu=\la u$ can be rewritten as
\[
\frac{\pd u_m}{\pd x^l}-\frac{\pd u_l}{\pd x^m}=\la\,(\mu_g)_{lms}u^s\,.
\]
Substituting this equation into~\eqref{DG*du} and utilizing the identity $(\mu_g)_{ijk}(\mu_g)_{lmn}g^{kn}=g_{il}g_{jm}-g_{im}g_{jl}$, the claim follows.
\end{proof}

In the proof of Theorem~\ref{T.trans} we will need the following lemma. Locally, it is an immediate consequence of the flow box theorem, while the global statement follows from a standard coordinate patching argument.

\begin{lemma}\label{L.dense}
  Let $w\in C^{r}(M,\La^1)$, $r\geq1$, and consider a compact subset $K\subset M\minus w^{-1}(0)$. Then
  \[
\big\{v\in C^{r}(M,\La^1):\supp v\subset K\big\}\subset \big\{Tw:T^i_j=g^{ik}h_{kj},\; h\in\cS^{r}(M)\big\}\,.
\]
\end{lemma}

For any $\bar g\in\cG^r(M)$, let us define the orthogonal projection $P^{\bar g}$ onto $\cK^{\perp_{\bar g}}$, which is a self-adjoint operator in  $(L^2(M,\La^1),\lan\cdot,\cdot\ran_{\bar g})$ that will be of considerable use in what follows. A convenient property of these projections is presented in the following

\begin{lemma}\label{L.P}
 For any $g,\bar g\in\cG^r(M)$, the bounded operator $P_g^{\bar g}:=P^{\bar g}|_{\cK^{\perp_g}}:\cK^{\perp_g}\to\cK^{\perp_{\bar g}}$ is bijective. Moreover, if $S\subset\cK^{\perp_g}$ is a closed subspace of finite codimension $n$, then $P_g^{\bar g}(S)$ is closed and has codimension $n$ in $\cK^{\perp_{\bar g}}$.
\end{lemma}
\begin{proof}
  To show that $P_g^{\bar g}$ is injective, notice that if $P_g^{\bar g}w=0$ for some $w\in\cK^{\perp_g}$, it then follows that $w\in\cK$, which implies that $w=0$. One can similarly establish that $P_g^{\bar g}$ is onto.

  Let us now suppose that we have the direct sum $\cK^{\perp_g}=S\oplus V$, with $\dim V=n$. Since $P_g^{\bar g}$ is onto, obviously $\cK^{\perp_{\bar g}}=P_g^{\bar g}(S)+P_g^{\bar g}(V)$, while the fact that the latter sum is direct is implied by the injectivity of $P_g^{\bar g}$. Since $P_g^{\bar g}(V)$ is closed and of dimension $n$, the statement follows.
\end{proof}

A key role in our analysis of the eigenvalues of the Beltrami operator will be played by the smooth map $\Phi: E×\RR\to L^2(M,\La^1)$ given by
\begin{equation}\label{Phi}
\Phi(g,u,\la):=(*_gd-\la)u
\end{equation}
and by the associated map $\Phi^{\bar g}:E×\RR\to \cK^{\perp_{\bar g}}$ defined as
\[
\Phi^{\bar g}(g,u,\la):=P^{\bar g}(*_gd-\la)u\,.
\]
Here $\bar g$ is any metric on $M$ of class $C^r$. To state the following lemma, let us recall that a linear map between Banach spaces is {\em Fredholm} if it has closed image and finite-dimensional kernel and cokernel.  A $C^1$ map between Banach manifolds $\Psi:N\to S$ is a {\em Fredholm map} if its derivative is Fredholm at every point of $N$. The {\em index} of a Fredholm operator is the difference between the dimensions of its kernel and cokernel.

\begin{lemma}\label{L.Fredholm}
  For each $g,\bar g\in \cG^r(M)$,
\[
\Phi^{\bar g}_g:=\Phi^{\bar g}(g,\cdot,\cdot):E_g×\RR\to\cK^{\perp_{\bar g}}
\]
is a Fredholm map of index $0$.
\end{lemma}
\begin{proof}
  We shall prove that $(D\Phi^{\bar g}_g)_{(u,\la)}:T_uE_g×\RR\to\cK^{\perp_{\bar g}}$ is a Fredholm operator for every $u\in E_g$ and $\la\in\RR$. To begin with, let us notice that
  \[
T_uE_g=\big\{v\in H^1(M,\La^1)\cap \cK^{\perp_g}: v\perp_g u\big\}
  \]
and  consider the linear map $F:T_uE_g×\RR\to \cK^{\perp_g}$ defined by
  \[
F(v,\nu):=(*_gd-\la)v-\nu u\,.
\]

Let us consider the self-adjoint operator $L$ on $(\cK^{\perp_g},\lan\cdot,\cdot\ran_g)$ given by 
\[
Lv:=(*_gd-\la)v\,,
\]
with domain $H^1(M,\La^1)\cap \cK^{\perp_g}$. We define $m:=\dim\ker L<\infty$, so that $m\neq0$ if and only if $\la$ is a nonzero eigenvalue of $*_gd$. If $u\perp_g\ker L$, then one can decompose $\cK^{\perp_g}$ into the $g$-orthogonal direct sum
\[
\cK^{\perp_g}=\ker L\oplus_g \RR u\oplus_g W\,,
\]
where $W$ is a closed subspace of codimension $m+1$. By the Fredholm alternative, it is standard that $L|_{(\RR u\,\oplus_g W)\cap H^1(M,\La^1)}$ is injective and has range $(\ker L)^{\perp_g}\cap \cK^{\perp_g}$ in this case. If $u\not\perp_g\ker L$, then
\begin{equation*}
\cK^{\perp_g}=N \oplus_g \RR u\oplus_g \widetilde W\,,
\end{equation*}
where $N\subset\ker L$ and $\widetilde W$ is a closed subspace of codimension $m$. By the Fredholm alternative, $L|_{\widetilde W\cap H^1(M,\La^1)}$ is injective and has range $(\ker L)^{\perp_g} \cap \cK^{\perp_g}$. Since $\im F$ is $L(W)+\RR u$ and $L(\widetilde W)+\RR u$ in each case, it is apparent that $F$ has closed image and finite dimensional kernel and cokernel. 

Notice that
\[
(D\Phi^{\bar g}_g)_{(u,\la)}(v,\nu)=P_g^{\bar g}F(v,\nu)\,.
\]
By Lemma~\ref{L.P} and the fact that $F(v,\nu)\in\cK$ if and only if $F(v,\nu)=0$, $\Phi_g^{\bar g}$ is a Fredholm map and $\ind(F)=\ind(\Phi_g^{\bar g})$ is independent of $u$ and $\la$.
In order to compute the index of $F$, we can assume that $u$ is not orthogonal to $\ker L$, so that the kernel of $F$ consists of the points $(v,\nu)\in T_uE_g×\RR$ such that
\begin{equation}\label{eq.Fredholm}
F(v,\nu)=Lv-\nu u=0\,.
\end{equation}
Since $u$ is not orthogonal to $\ker L$, this implies that $\nu=0$, so that $\ker F=N×\{0\}$ has dimension $m-1$. A similar argument shows that $\im F=L(\widetilde W)\oplus \RR u=(\ker L)^{\perp_g}\cap \cK^{\perp_g}\oplus\RR u$ has codimension $m-1$ in $\cK^{\perp_g}$, which yields that $\ind(F)=0$, as claimed.
\end{proof}

\section{Analysis of the Beltrami operator}
\label{S.Beltrami}

In this section we aim to prove that, for a residual set of $C^r$ metrics, all the eigenvalues of $*_gd$ are simple and the zeros of the corresponding eigenfunctions are hyperbolic. As in~\cite{Uh76}, our approach relies on techniques from infinite dimensional transversality theory, but the peculiarities of the Beltrami equation induce various fundamental complications that were totally absent in Uhlenbeck's analysis of scalar second order elliptic equations. Let us begin by discussing the new technicalities in some detail using the notation of the previous section without further mention.

Infinite dimensional transversality theory is tailored for Fredholm maps, whose kernel and cokernel are of finite dimension, since this feature allows one to invoke the Sard--Smale theorem. Here we need to apply these ideas to the map $\Phi$, which is the central object in our proof of the generic simplicity of the spectrum of the Beltrami operator. The obvious fact that the symbol of the Beltrami operator is not elliptic manifests itself in the necessity of defining $\Phi$ on a geometrically nontrivial fiber bundle instead of on a Cartesian product of Banach manifolds (cf.~\eqref{Phi}), which makes the analysis of the differential of $\Phi$ nontrivial due to the constraints that determine the tangent space~\eqref{TE}. A closely related problem is that the image of the differential of $\Phi_g$ has infinite codimension. Actually, $\Phi$ does not have any nontrivial regular values, and one can only prove that it is transverse to the infinite dimensional submanifold $\cK$ (Corollary~\ref{C.trans}). Exploiting this fact using a parametric transversality argument is usually rather problematic, although in this case we manage to get by using Lemmas~\ref{L.P} and~\ref{L.Fredholm}.

The vectorial character of the Beltrami equation makes the study of the generic nondegeneracy of the zero set of the eigenfunctions rather different than in the case of scalar second order equations. Indeed, this problem is obviously equivalent to establishing the surjectivity of the differential of an evaluation map, which in the scalar case can be accomplished through a simple rescaling argument~\cite{Uh76}. As we shall see in Theorem~\ref{T.sing}, the argument for the Beltrami operator is considerably more involved.

It is well known that an eigenfunction of the Beltrami operator satisfies the incompressible steady Euler equations, and that it defines a contact structure on the manifold provided that its zero set is empty. Consequently, the eigenfunctions of the Beltrami operator have been thoroughly studied in the context of hydrodynamics by several authors, and more recently in connection with contact geometry by Etnyre and Ghrist~\cite{EG00}. It should be mentioned that, as a technical tool for proving generic
hydrodynamic instability, the latter authors have
claimed~\cite[Theorem 2.1 and Lemma 2.2]{EG05} that for a generic
$C^r$ metric the eigenvalues of the Beltrami operator are simple and
the eigenfunctions have nondegenerate zeros. Unfortunately, in this interesting paper the authors provide a proof of this fact which contains some substantial gaps. Being of direct
interest for this section, we shall next outline the main points that remain to be
fixed in~\cite{EG05}; a complete, self-contained proof will be given in Sections~\ref{SS.spec} and~\ref{SS.eigenfunctions}.

In~\cite{EG05}, the authors essentially consider the function
\[
\widetilde\Phi(g,u,\la):=\big( g,(*_gd-\la)u\big)
\]
mapping $E×\RR$ into the vector bundle $E'$ over $\cG^r(M)$ with fiber $\cK^{\perp_g}$, and claim that $\widetilde\Phi$ is transverse to the zero section of $E'$. A quick inspection reveals that their proof of this result is not conclusive: they do not address the surjectivity of $(D\widetilde\Phi)_{(g,u,\la)}$, as required, but only the density of its image. Moreover, they do not show why the image of $(D\widetilde\Phi)_{(g,u,\la)}$ is dense, which at least requires to analyze the variation of the Beltrami operator with the metric using the description of $T_{(g,u)}E$. These computations are important, as they are precisely the reason why Uhlenbeck's results have not been extended to the Hodge Laplacian (and cannot be, at least in even dimension and middle degree~\cite{Mi80}).

A second mistake, which is the crucial one and directly related to the first, is that their proof of the nondegeneracy of the zeros of the eigenfunctions hinges on their claim that
\[
\big\{D(*d)_g(h)u: h\in\cS^r(M)\big\}\supset C^r(M,\La^1)
\]
for any eigenfunction $u$. This is obviously not true because the left hand side vanishes on $u^{-1}(0)$, as revealed by the computations leading to Lemma~\eqref{L.D}. 

A third gap in the proof is that Uhlenbeck's parametric transversality Theorem~1~\cite{Uh76} does not apply to $\widetilde\Phi$ because it is transverse to the infinite-dimensional submanifold $\cG^r(M)×\{0\}$, not to a point. Our treatment of this issue makes use of some other results  we prove in this section. It should be noticed that the ideas developed in this article could be used to derive a proof (of similar complexity) of the generic simplicity of the spectrum and nondegeneracy of the zeros of the eigenfunctions of the Beltrami operator using $\widetilde\Phi$.

\subsection{Spectrum of the Beltrami operator}
\label{SS.spec}

Let us now move on to the main technical result needed to establish the generic simplicity of the eigenvalues of the Beltrami operator. Before stating this theorem, let us recall that a point $x\in N$ is a {\em regular point} of a $C^1$ map $\Psi:N\to S$ between Banach manifolds if the derivative $D\Psi_x:T_xN\to T_{\Psi(x)} S$ is onto. A point $y\in S$ is a {\em regular value} of $\Psi$ if $\Psi^{-1}(y)=\emptyset$ or all the points in its preimage are regular. More generally, the map $\Psi$ is said to be {\em transverse} to a submanifold $\Si\subset S$ if
\[
{\rm im}(D\Psi_x) + T_{\Psi(x)}\Si = T_{\Psi(x)}S
\]
for all $x\in N$ such that $\Psi(x)\in\Si$.

\begin{theorem}\label{T.trans}
For each $\bar g\in\cG^r(M)$, $0$ is a regular value of $\Phi^{\bar g}$.
\end{theorem}
\begin{proof}
Let us take $(g,u,\la)\in (\Phi^{\bar g})^{-1}(0)$. By the definition of~\eqref{Eg} and standard regularity results, this actually implies that $u$ is an eigenfunction of $*_gd$ with eigenvalue $\la\neq0$ and of class  $C^{r,\al}$ for all $\al<1$.

Let us take two $1$-forms $\al,\be\in H^1(M,\La^1)$ and consider their scalar product
\[
\lan\al,\be\ran_g=\int g^{ij}\,\al_i\,\be_j\,d\mu_g
\]
as a function of $g$. In local coordinates, the volume element reads as $d\mu_g=|g|^{1/2}dx^1\,dx^2\,dx^3$, so a straightforward computation as in Lemma~\ref{L.D} shows that its derivative $(D\lan\al,\be\ran)_g:\cS^r(M)\to\RR$ is
\begin{equation}\label{Dlan}
(D\lan\al,\be\ran)_g(h)=\int\bigg(\frac{\tr_gh}2g^{ij}-h^{ij}\bigg)\,\al_i\,\be_j\, d\mu_g\,.
\end{equation}
Since
\[
E=\big\{(g,u): \|u\|_g=1,\;\lan u,\al\ran_g=0\;\text{for all }\al\in \cK\big\}\,,
\]
with $(g,u)\in\cG^r(M)× H^1(M,\La^1)$, it follows from~\eqref{Dlan} that
\begin{align}
T_{(g,u)}E=\bigg\{(h,v): \;&\lan\al,v\ran_g+ \notag \int\bigg(\frac{\tr_gh}2g^{ij}-h^{ij}\bigg)\,\al_i\,u_j\, d\mu_g=0\;\text{for all }\al\in \cK\,,\\
&2\lan u,v\ran_g+ \int\bigg(\frac{\tr_gh}2g^{ij}-h^{ij}\bigg)\,u_i\,u_j\, d\mu_g=0\bigg\}\,,\label{TE}
\end{align}
where $(h,v)\in\cS^r(M)× H^1(M,\La^1)$. Let us consider $\ga\in L^2(M,\La^1)$ orthogonal to $\im D\Phi_{(g,u,\la)}$, i.e., such that
\begin{equation}\label{orth}
\big\lan (*_gd-\la)v-\nu u+D(*d)_g(h)u,\ga\big\ran_g=0
\end{equation}
for all $(h,v)\in T_{(g,u)}E$, $\nu\in\RR$. Let us decompose $\ga=\al+\be$, with $\al\in\cK$ and $\be\in\cK^{\perp_g}$. Taking $h=0$ and $\nu=0$ in the latter equation, we immediately derive that
\[
\lan(*_gd-\la)v, \ga\ran_g=0
\]
for all $v\in H^1(M,\La^1)\cap \cK^{\perp_g}$, which implies that
\[
\lan(*_gd-\la)v, \be\ran_g=0
\]
for all $v\in H^1(M,\La^1)$. It then follows that $\be\in C^{r,\al}(M,\La^1)$ solves the equation
\begin{equation}\label{orth2}
(*_gd-\la)\be=0\,.
\end{equation}

Let us now consider arbitrary variations in~\eqref{orth} with $\nu=0$, so that by Eq.~\eqref{orth2} we have
\begin{align}\notag
0&=\big\lan\al+\be,(*_gd-\la)v+D(*d)_g(h)u\big\ran_g\\
&=\big\lan\al, D(*d)_g(h)u-\la v\big\ran_g +\big\lan\be, D(*d)_g(h)u\big\ran_g\,.\label{zero}
\end{align}
Since $\al\in\cK$, it follows from~\eqref{TE} and Lemma~\ref{L.D} that
\begin{align*}
\big\lan D(*d)_g(h)u-\la v,\al\big\ran_g&= \big\lan D(*d)_g(h)u,\al\big\ran_g +\la\int\bigg(\frac{\tr_gh}2g^{ij}-h^{ij}\bigg)\,\al_i\,u_j\, d\mu_g=0\,,
\end{align*}
so that~\eqref{zero} yields
\begin{equation}\label{orth3}
\int \bigg(h^{ij}-\frac{\tr_gh}2g^{ij}\bigg)\be_iu_j\,d\mu_g=0
\end{equation}
for all $h\in \cS^r(M)$.

For any symmetric tensor $T\in\cS^r(M)$, we infer that
\begin{equation}\label{orth4}
\lan Tu,\be\ran_g=0
\end{equation}
by taking $h=T-(\tr_gT)g$ in Eq.~\eqref{orth3}. Since $(*_gd-\la)u=0$, it is a trivial observation that the co-closed 1-form $u$ is an eigenfunction of the Hodge Laplacian~$\Delta_g$ with eigenvalue $\la^2$, and therefore it does not vanish in any (nonempty) open subset of $M$ by the unique continuation theorem~\cite{Ar57,Ka88}. From Lemma~\ref{L.dense} it then follows that $\{Tu:T\in\cS^r(M)\}$ is dense in $L^2(M,\La^1)$, so that Eq.~\eqref{orth4} implies that $\be=0$.

The above argument shows that $P^g(\im D\Phi_{(g,u,\la)})$ is dense in $\cK^{\perp_g}$. Since $\Phi_g^g$ is Fredholm by Lemma~\ref{L.Fredholm}, $P^g(\im D(\Phi_g)_{(u,\la)})$ has finite codimension $n$ in $\cK^{\perp_g}$, and therefore we can take an $n$-dimensional subspace $V\subset \cS^r(M)$ such that
\[
\cK^{\perp_g}=P^g\big(\im D(\Phi_g)_{(u,\la)}\oplus  D\Phi_{(g,u,\la)}(V)\big)\,.
\]
Hence we deduce that $P^g(\im D\Phi_{(g,u,\la)})=\cK^{\perp_g}$, which implies that $D(\Phi^{\bar g})_{(g,u,\la)}$ is onto since
\[
\im D(\Phi^{\bar g})_{(g,u,\la)}=P^{\bar g}\big(\im D\Phi_{(g,u,\la)}\big)= P^{\bar g}_g(\cK^{\perp_g})=\cK^{\perp_{\bar g}}\,.
\]
Lemma~\ref{L.P} has been used to derive the final identity. By definition, $0$ is then a regular value of $\Phi^{\bar g}$, completing the proof of the theorem.
\end{proof}

\begin{remark}
  As noted by Uhlenbeck in the case of the scalar Laplacian~\cite{Uh76}, one can replace the set of $C^r$ metrics by the set of $C^r$ metrics which differ from a fixed one only on a given (nonempty) open subset of $M$. Indeed, one can easily check that the proofs of our genericity results, including Theorem~\ref{T.main}, remain valid when the metrics we consider are assumed to coincide with a fixed metric $\bar g$ but in a (nonempty) open set. What should be noticed here is that, contrary to what happens with the scalar Laplacian, in the case of differential forms it is not sufficient to consider global conformal deformations of the metric to derive the desired results.
\end{remark}

As an aside, observe that from the proof of the theorem one also obtains the following 

\begin{corollary}\label{C.trans}
  $\Phi$ is transverse to $\cK$.
\end{corollary}

We are now ready to establish the generic simplicity of the eigenvalues of the Beltrami operator making use of the transversality theorem we present below. In order to state it, let us recall that a subset of a topological space is {\em residual} if it is the countable intersection of open dense sets. It is standard that a residual set of a complete metric space is necessarily dense by the Baire category theorem.

\begin{theorem}[\cite{Qu70,Uh76}]\label{T.Uhl}
Let $\pi:N\to B$ be a fiber bundle and let $\Psi:N\to S$ be a $C^{1}$ map, with $N,B,S$ Banach manifolds and $B,S$ separable. If $y$ is a regular value of $\Psi$ and $\Psi|_{\pi^{-1}(b)}$ is a Fredholm map of index $0$ for all $b\in B$, then the set
  \[
\big\{b\in B:y\text{ is a regular value of }\Psi|_{\pi^{-1}(b)}\big\}
  \]
  is residual in $B$. 
\end{theorem}

\begin{theorem}\label{T.values}
  The eigenvalues of the Beltrami operator $*_gd$, acting on $H^1(M,\La^1)\cap\cK^{\perp_g}$, are all simple for a residual set of $C^r$ metrics.
\end{theorem}
\begin{proof}
By Theorem~\ref{T.trans}, $0$ is a regular value of $\Phi^{\bar g}$, so that Theorem~\ref{T.Uhl} ensures that $0$ is also a regular value of $\Phi_g^{\bar g}$ for a residual subset of $\cG^r(M)$. Since $\Phi_g^g=P^g_g\Phi_g^{\bar g}$, Lemma~\ref{L.P} ensures that $(D\Phi_g^{\bar g})_{(u,\la)}$ is onto if and only if $(D\Phi_g^{g})_{(u,\la)}$ is. 

Let us fix a metric $g$ in this residual subset and suppose that $u,\hat u\in E_g$ are two linearly independent eigenfunctions of $*_gd$ with eigenvalue $\la$. As $D(\Phi_g^g)_{(u,\la)}:T_uE_g×\RR\to \cK^{\perp_g}$ is onto, there must exist $(v,\nu)\in T_uE_g×\RR$ such that
\[
D(\Phi_g^g)_{(u,\la)}(v,\nu)=(*_gd-\la)v-\nu u=\hat u\,.
\]
Since $\hat u\in\ker(*_gd-\la)$ is not proportional to $u$, this is absurd.
\end{proof}

\subsection{Eigenfunctions of the Beltrami operator}
\label{SS.eigenfunctions}

We shall next show that the zero set of the eigenfunctions of the Beltrami operator generically consists of hyperbolic singular points.  The proof of this fact utilizes the following transversality theorem. Before presenting the statement of these results, let us recall that, given a $1$-form $u\in C^2(M,\La^1)$, a singular point $x\in u^{-1}(0)$ is {\em hyperbolic} if all the eigenvalues of the Jacobian matrix $\nabla u(x)$ of $u$ at $x$ have nonzero real part. In particular, the latter condition implies that $x$ is an isolated zero of $u$.

\begin{theorem}[\cite{Uh76}]\label{T.Uhl2}
Let $Z,B,N_1,N_2,N_2'$ be separable Banach manifolds, with $N_2'\subset N_2$ and $N_1,N_2$ of finite dimension. Let $\pi:Z\to B$ be a $C^k$ Fredholm map of index $0$. If $\si:Z× N_1\to N_2$ is of class $C^k$ and transverse to $N_2'$, with $k>\max\{0,\dim N_1+\dim N_2'-\dim N_2\}$, then the set 
\[
\big\{ b\in B: \si|_{\pi^{-1}(b)} \text{ is transverse to } N_2'\big\}
\]
is residual in $B$.
\end{theorem}

\begin{theorem}\label{T.sing}
  The zeros of the eigenfunctions of the Beltrami operator $*_gd$, acting on $E_g$, are all hyperbolic for a residual set of $C^r$ metrics.
\end{theorem}
\begin{proof}
  By Corollary~\ref{C.trans}, $Z:=\Phi^{-1}(\cK)$ is a smooth submanifold of $E× \RR$. Since clearly $\Phi^{-1}(\cK)=\Phi^{-1}(0)$ because $\Phi(\{g\}× E_g×\RR)\perp_g \cK$, it follows that $(g,u,\la)$ lies in $Z$ if and only if $u$ is an eigenfunction of $*_gd$ with nonzero eigenvalue $\la$ and unit norm.

  Let us also consider the natural projection $\Pi: Z\to\cG^r(M)$. It is not difficult to see that $\Pi$ is a Fredholm map of index $0$. Indeed, let us take $(g,u,\la)\in Z$, so that $D\Pi_{(g,u,\la)}: T_{(g,u,\la)}Z\to \cS^r(M)$ is given by
  \begin{equation}\label{DPi}
D\Pi_{(g,u,\la)}(h,v,\nu)=h\,,
\end{equation}
with
  \begin{align}\notag
    T_{(g,u,\la)}Z:=\bigg\{(h,v,\nu):\; &(*_gd-\la)v -\nu u+D(*d)_g(h)u=0,\\
    &2\lan u,v\ran_g+ \int\bigg(\frac{\tr_gh}2g^{ij}-h^{ij}\bigg)\,u_i\,u_j\, d\mu_g=0\bigg\}\label{TZ}
  \end{align}
and $(h,v,\nu)\in \cS^r(M)× H^1(M,\La^1)× \RR$. This expression should be compared with~\eqref{TE}: here, we do not need to include the constraint associated to the orthogonality of $u$ and $\cK$ because it is automatically implied by the fact that $u\in Z$.
  
It is apparent from Eq.~\eqref{DPi} that the image of $D\Pi_{(g,u,\la)}$ is closed.
Since $\Phi$ is transverse to $\cK$ by Corollary~\ref{C.trans}, it is not difficult to derive (cf.\ e.g.~\cite{Qu70}) that
\begin{align}\notag
  \codim_{\cS^r(M)}\im D\Pi_{(g,u,\la)} &= \codim_{T_{(g,u,\la)}(E×\RR)} \big( T_{(u,\la)}(\{g\}× E_g×\RR)+ T_{(g,u,\la)}Z\big)\\
    & =\codim_{L^2(M,\La^1)}\big(\im(D\Phi(g,\cdot,\cdot))_{(u,\la)}+\cK\big)\notag\\
    &=\codim_{\cK^{\perp_{\bar g}}}\im (D\Phi_g^{\bar g})_{(u,\la)}\,,\label{eq.Fred2}
\end{align}
where $\bar g$ is an arbitrary metric in $\cG^r(M)$ and $\Phi_g^{\bar g}$ is defined as in Lemma~\ref{L.Fredholm}.

By~\eqref{DPi} and~\eqref{TZ}, the kernel of $D\Pi_{(g,u,\la)}$ consists of the points $(0,v,\nu)$ which satisfy the equations
\[
(*_gd-\la)v=\nu u\,,\qquad v\perp_g u\,.
\]
From Eq.~\eqref{eq.Fredholm} it then follows that
\[
\ker D\Pi_{(g,u,\la)}=\ker (D\Phi_g^{\bar g})_{(u,\la)}\,,
\]
so that, by~\eqref{eq.Fred2} and Lemma~\ref{L.Fredholm},
\[
\ind\Pi=\ind\Phi_g^{\bar g}=0\,,
\]
as we wanted to show.

If $w\in H^2(M,\La^*):=\bigoplus_{p=0}^3 H^2(M,\La^p)$ is a differential form $g$-orthogonal to $\ker d$, the action of the elliptic differential operator $d+\de_{g}$ satisfies
\[
\big\|(d+\de_{g})^jw\big\|_{g}=\big\|(*_{g}d)^jw\big\|_{g}
\]
for $j\leq2$. Therefore, it is an easy observation that any $H^2$ norm on $H^2(M,\La^1)\cap \cK^{\perp_g}$ is equivalent to the norm
\[
|w|_{2,g}:=\sum_{j=0}^2\big\|(*_gd)^jw\big\|_g\,.
\]
Since $*_gd: H^1(M,\La^1)\cap\cK^{\perp_g}\to\cK^{\perp_g}$ is self-adjoint, there exists an orthonormal basis of $C^{r,\al}$ eigenfunctions $(\vp_n)_{n=1}^\infty\subset \cK^{\perp_g}$ with (not necessarily distinct) eigenvalues $\la_n$. It then follows from the previous argument that, for any $w\in H^{2}(M,\La^1)\cap \cK^{\perp_g}$, the eigenfunction expansion
\begin{equation}\label{expansion}
\sum_{n=1}^N\lan w,\vp_n\ran_g\,\vp_n
\end{equation}
converges to $w$ in $H^{2}(M,\La^1)$ as $N\to\infty$. In particular, this eigenfunction expansion converges pointwise by the Sobolev embedding theorem.

Let us now consider the evaluation map $\ev: Z× M\to T^*M$, defined by
\[
\ev(g,u,\la,x):=u(x)\,.
\]
The $C^{r,\al}$ regularity of the eigenfunctions readily imply that
$\ev$ is of class $C^{r}$. Our next goal will be to prove that $\ev$
is transverse to the zero section of $T^*M$. To this end, let us take
any $(g,u,\la,x)\in Z× M$ such that $u(x)=0$ and show that
$D\ev_{(g,u,\la,x)}$ is transverse to the zero section. This is
tantamount to showing that $0$ is a regular value of the map
$\ev_x:Z\to T^*_xM$ given by
\[
\ev_x(g,u,\la):=u(x)\,.
\]

As
\[
(D\ev_x)_{(g,u,\la)}(h,v,\nu):=v(x)\,,
\]
the transversality of $\ev$ to the zero section will follow once we prove that
\[
V_x:=\big\{v(x): (h,v,\nu)\in T_{(g,u,\la)}Z\big\}
\]
is actually the whole space $T^*_xM$. By linearity, $V_x$ is a linear subspace of $T^*_xM$, so a necessary and sufficient condition for $V_x$ not to be equal to $T_x^*M$ is that there exists a nonzero $\xi\in T_x^*M$ such that the inner product
\begin{equation}\label{xi}
\xi\cdot \eta:=g^{ij}(x)\,\xi_i\eta_j
\end{equation}
is zero for all $\eta\in V_x$. 

Let us now introduce the resolvent operator
\[
R_\la w:=\sum_{\la_n\neq\la}\frac{\lan w,\vp_n\ran_g}{\la_n-\la}\,\vp_n\,,
\]
which satisfies
\[
(*_gd-\la)R_\la w=w
\]
for any $w\perp_g\ker(*_gd-\la)$, $w\in\cK^{\perp_g}$. As argued in the proof of Theorem~\ref{T.trans}, Lemmas~\ref{L.D} and~\ref{L.dense} imply that for any $f\in C^r(M,\La^1)$ such that $\supp(f)\cap u^{-1}(0)=\emptyset$ there exists a symmetric tensor $h_f\in\cS^r(M)$ such that
\[
f=-D(*d)_g(h_f)u\,.
\]
Let us take $f$ orthogonal to $\ker(*_gd-\la)$ and $\cK$. In this case, the characterization of the tangent space~\eqref{TZ} and the fact that $u$ belongs to the kernel of $*_gd-\la$ ensure that the element $(h_f,v_f,0)$ belongs to $T_{(g,u,\la)}Z$, where
  \[
v_f:=R_\la f-c_f u
  \]
  and
  \[
c_f:=\frac12\int\bigg(\frac{\tr_gh_f}2g^{ij}-(h_f)^{ij}\bigg)\,u_i\,u_j\, d\mu_g\,.
\]
Moreover, the partial sums
\begin{equation}\label{expansion2}
\sum_{n\leq N,\,\la_n\neq\la}\frac{\lan f,\vp_n\ran_g}{\la_n-\la}\,\vp_n(x)
\end{equation}
converge to $v_f(x)+c_fu(x)$ pointwise as $N\to\infty$ by the argument used to prove the pointwise convergence of~\eqref{expansion} and the fact that the eigenvalues of $*_gd$ do not accumulate at $\la$.

Since $u(x)=0$, one obviously has that $v_f(x)=R_\la f(x)$. Eq.~\eqref{xi} and the pointwise convergence of~\eqref{expansion2} then imply that
\[
\xi\cdot v_f(x)=\sum_{\la_n\neq\la} \frac{\lan\vp_n,f\ran_g}{\la_n-\la}\xi\cdot \vp_n(x)=0
\]
for all $f$ as above. By the unique continuation theorem~\cite{Ar57,Ka88}, $u^{-1}(0)$ is a closed set of empty interior, so that the set of admissible $f$ is dense in $\ker(*_gd-\la)^{\perp_g}\cap \cK^{\perp_g}$. It then follows that
\begin{equation}\label{perpxi}
\xi\cdot \vp_n(x)=0
\end{equation}
whenever $\la_n\neq\la$. As $\xi\cdot u(x)=0$ is automatic and $v(x)\in V_x$ if $v\in\ker(*_gd-\la)$ and $v\perp_g u$ by Eq.~\eqref{TZ}, it stems that~\eqref{perpxi} holds true for all $n\in\NN$. Hence, the pointwise convergence of~\eqref{expansion} and Eq.~\eqref{perpxi} now imply that
\[
\xi\cdot w(x)=0
\]
for all $w\in H^2(M)\cap \cK^{\perp_g}$, which is absurd.

Therefore we infer that the evaluation map is transverse to the zero section of $T^*M$.
We can now apply Theorem~\ref{T.Uhl2} to the maps $\Pi:Z\to \cG^r(M)$ and $\ev:Z× M\to T^*M$ to show that the set
\[
\big\{ g\in \cG^r(M): \ev|_{\Pi^{-1}(g)} \text{ is transverse to the zero section of } T^*M\big\}
\]
is residual in $\cG^r(M)$. Notice that the above transversality condition simply means that all the zeros of the eigenfunctions of $*_gd$ with nonzero eigenvalue are nondegenerate. Since $\de_g u=0$, a straightforward computation shows that the real part of the eigenvalues of the matrix $\nabla u(x)$ must be nonzero as well, for any $x\in u^{-1}(0)$.
\end{proof}

\section{The Hodge Laplacian}
\label{S.Hodge}

In this section we shall utilize the results on the Beltrami operator derived in the previous sections to prove the main theorem. In doing this, we make essential use of the fact that $M$ is three-dimensional.

It should be remarked that we have chosen to base our approach on an analysis of the Beltrami operator because highly nontrivial complications arise when one tries to directly apply Uhlenbeck's method to the Laplacian on $1$-forms. Indeed, the action of $\De$ on co-exact forms is given by $(*_gd)^2$, and if  $u_±\in E_g$ are eigenfunctions of $*_gd$ with nonzero eigenvalue $±\la_0$, it follows that the map $\widehat\Phi:E×\RR\to L^2(M,\La^1)$ defined by $\widehat\Phi(g,u,\la):=(\De_g-\la)u$ cannot be transverse to $\cK$ because
\begin{align*}
\lan u_-,(D\widehat\Phi)_{(g,u_+,\la_0^2)}(h,v,\nu)\ran_g&=\lan u_-, (\De_g-\la_0^2)v-\nu u_++ D(\De)_g(h)u_+\ran_g\\
&=\lan u_-,D((*d)^2)_g(h)u_+\ran_g\\
&=\lan*_gd u_-,D(*d)_g(h)u_+\big\ran_g +\big\lan u_-,D(*d)_g(h)*_gdu_+\ran_g\\
&=0\,.
\end{align*}

For the case of the scalar Laplacian, a perturbation-theoretic approach can also be used to prove the genericity of simple eigenvalues~\cite{Al78,BU83}. Nonetheless, as one can infer from~\cite{CPR01}, it is not obvious at all how this technique can be adapted to the case of $1$-forms, as conformal variations of the metric are a priori not sufficient to break the degeneracy of eigenvalues.

\begin{proof}[Proof of Theorem~\ref{T.main}]
Let us suppose that $\la$ is a nonzero simple eigenvalue of the self-adjoint operator $*_gd: H^1(M,\La^1)\cap \cK^{\perp_g}\to\cK^{\perp_g}$, and let $u$ be an associated eigenfunction of norm $1$. In this case, it is then standard~\cite{Ka80} that the eigenvalues and eigenfunctions depend smoothly on the metric in a neighborhood of $(g,u,\la)$. More precisely, there exists a neighborhood $N_{g,u,\la}\subset\cG^r(M)$ and smooth functions $\ell: N_{g,u,\la}\to\RR$ and $U: N_{g,u,\la}\to H^1(M,\La^1)$ such that
\begin{equation}\label{cogn}
*_{\bar g}dU(\bar g)=\ell(\bar g)\,U(\bar g)
\end{equation}
for all $\bar g\in N_{g,u,\la}$. Moreover, $\ell(g)=\la$, $U(g)=u$ and $U(\bar g)\in E_{\bar g}$ for all $\bar g\in N_{g,u,\la}$.
Taking derivatives in~\eqref{cogn} with respect to the metric and evaluating at $g$, we arrive at the equation
\[
D(*d)_g(h)u+(*_gd-\la)(DU)_g(h)-(D\ell)_g(h)u=0\,,
\]
which readily yields
\begin{equation}\label{cogn2}
(D\ell)_g(h)=\la\int\bigg(h^{ij}-\frac{\tr_g h}2g^{ij}\bigg)\, u_i\, u_j\,d\mu_g
\end{equation}
by taking the inner product with $u$ and using Lemma~\ref{L.D}.

Similarly, if $\si$ is a nonzero simple eigenvalue of the scalar Laplacian $\De_g:=\de_gd: H^2(M)\to L^2(M)$ with normalized eigenfunction $f$, there exists a neighborhood $R_{g,f,\si}\subset\cG^r(M)$ and smooth functions $s: R_{g,f,\si}\to\RR$ and $F: R_{g,f,\si}\to H^2(M)$ such that
\begin{equation*}
\De_{\bar g}F(\bar g)=s(\bar g)\,F(\bar g)\,,
\end{equation*}
$s(g)=\si$, $F(g)=f$ and $\int F(\bar g)^2\,d\mu_{\bar g}=1$. Proceeding as above, one immediately arrives at the formula
\begin{equation}\label{cogn3}
(Ds)_g(h)=-\int\bigg(\frac{\De\tr_g h}4f^2+h(\nabla f,\nabla f)\bigg)\,d\mu_g
\end{equation}
for the first variation of the eigenvalue.

For concreteness, let us label the nonzero eigenvalues $(\la_n(g))_{n=1}^\infty$ and $(\si_n(g))_{n=1}^\infty$ of the Beltrami operator $*_gd$ and the scalar Laplacian so that
\[
\la_{n+1}^2(g)\geq \la_{n}^2(g)\qquad\text{and}\qquad \si_{n+1}(g)\geq \si_n(g)\,.
\]
The associated normalized eigenfunctions will be denoted by $(u_n(g))_{n=1}^\infty$ and $(f_n(g))_{n=1}^\infty$. For each $k\in\NN$, let us define the sets of metrics
\begin{align*}
  \Ga_k^1&:=\big\{ g\in\cG^r(M): \la_n(g)\neq \la_m(g)\text{ for all } 1\leq n\neq m\leq k \big\}\,,\\
  \Ga_k^2&:=\big\{ g\in\cG^r(M): \si_n(g)\neq \si_m(g)\text{ and the zeros of } df_n(g) \text{ are nondegenerate}\\
  &\hspace{23em} \text{for all } 1\leq n\neq m\leq k \big\}\,,\\
  \Ga_k^3&:=\big\{ g\in\cG^r(M): \text{the zeros of } u_n(g) \text{ are hyperbolic for all } 1\leq n\leq k \big\}\,,\\
  \Ga_k^4&:=\big\{ g\in\cG^r(M): \la_n(g)\neq -\la_m(g)\text{ for all } 1\leq n,m\leq k \big\}\,,\\
  \Ga_k^5&:=\big\{ g\in\cG^r(M): \si_n(g)\neq \la_m^2(g)\text{ for all } 1\leq n,m\leq k \big\}\,.
\end{align*}

Our purpose is to show that each set $\Ga_k^a$ is open and dense. The fact that $\Ga_k^1$ is dense is an immediate implication of Theorem~\ref{T.values}, while the openness readily follows from Eq.~\eqref{cogn2} and the fact that the condition $\la_n(g)\neq\la_m(g)$ is stable. By Uhlenbeck's theorem~\cite{Uh76}, the set $\Ga_k^2$ is also dense, while its openness is apparent from Eq.~\eqref{cogn3} and the fact that its defining conditions are stable.
Similarly, the density and openness of $\Ga_k^3$ is immediate from Theorem~\ref{T.sing} and the $C^1$-stability of the hyperbolic zeros of a differential form.

Arguing as above, it is obvious that $\Ga_k^1\cap\Ga_k^4$ is open. In order to see that it is dense, let us assume that $g\in\Ga_k^1$ is such that $\la_{n_j}(g)=-\la_{m_j}(g)$, for some $(n_j)_{j=1}^N,\, (m_j)_{j=1}^N\subset\{1,\dots,k\}$. For each $T\in\cS^r(M)$, let us set
\[
H(T):=T-(\tr_gT)g\,.
\]
By~\eqref{cogn2},
\[
(D\ell_n)_g(H(T))=\la\int T^{ij}\, u_n(g)_i\, u_n(g)_j\,d\mu_g\,,
\]
where the function $\ell_n$ describes the variation of the $n$-th eigenvalue of $*_gd$ with respect to the metric. Defining
\[
\rho_j(T):=(D\ell_{n_j})_g(H(T))+(D\ell_{m_j})_g(H(T))\,,
\]
it stems from~\eqref{cogn2} and the fact that $u_{n_j}(g)$ and $u_{m_j}(g)$ are linearly independent functions in $L^2(M,\La^1)$ that there exists $T_j\in\cS^r(M)$ satisfying $\rho_j(T_j)\neq0$. Let us now take real numbers $c_j$ such that $c_1:=1$,
\[
|c_j|<\min_{1\leq l<j}\bigg|\frac{\rho_l(c_1T_1+\cdots+c_{j-1}T_{j-1})}{\rho_l(T_j)}\bigg|\quad \text{for } j\geq2
\]
and $c_j=0$ if and only if $\rho_j(c_1T_1+\cdots+c_{j-1}T_{j-1})\neq0$. By construction,
\[
\rho_j(c_1T_1+\cdots+c_{N}T_{N})\neq0
\]
for all $1\leq j\leq N$, which implies that there exists a metric $\bar g\in\Ga_k^4$ arbitrarily close to $g$ in the $C^r$ topology. The density of $\Ga_k^1\cap\Ga_k^4$ is then a straightforward consequence of the density of $\Ga_k^1$.

Let us now consider the set
\[
\Ga_k:=\bigcap_{a=1}^5\Ga_k^a\,,
\]
which again is clearly open. To prove that it is also dense, let us take $g\in\bigcap_{a=1}^4\Ga_k^a$ and suppose that $\si_n(g)=\la_m^2(g)$. If $\tr_gh=0$, the variations of $\si_n$ and $\la_n^2$ with the metric along $h$ are respectively given by
\begin{align}
  (D\ell^2_m)_g(h)&=2\la_m^2(g)\int h^{ij}\, u_m(g)_i\, u_m(g)_j\,d\mu_g\,,\label{cogn4}\\
  (Ds_n)_g(h)&=-\int h\big(\nabla f_n(g),\nabla f_n(g)\big)\,d\mu_g\label{cogn5}
\end{align}
by Eqs.~\eqref{cogn2} and~\eqref{cogn3}. Since the traceless tensor $h$ can be chosen to be zero but in an open subset where
\[
\big(h^{ij}u_m(g)_iu_m(g)_j\big)\,h\big(\nabla f_n(g),\nabla f_n(g)\big)>0\,,
\]
Eqs.~\eqref{cogn4} and~\eqref{cogn5} show that there exists an arbitrarily $C^r$-small deformation $\bar g$ of $g$ for which $\si_n(\bar g)\neq \la_m^2(\bar g)$. Arguing as in the previous paragraph, this can be readily seen to imply that $\Ga_k$ is also dense, as we wanted to prove.

By the Hodge decomposition and the fact that
\[
\De_gdf=d\De_g f\,,\qquad \De_g u=(*_gd)^2u
\]
for any exact $1$-form $df$ and co-exact $1$-form $u$, the nonzero eigenvalues of the Laplacian on $1$-forms are given by the union of the nonzero eigenvalues of the scalar Laplacian and the squared eigenvalues of the Beltrami operator $*_gd$, counting multiplicities. Therefore, the desired set of metrics can be taken to be
\[
\Ga:=\bigcap_{k=1}^\infty\Ga_k\,,
\]
which is residual because each $\Ga_k$ is open and dense. Indeed, the definitions of $\Ga$ and $\Ga_k$ imply that the eigenvalues and eigenfunctions of the Laplacian on $0$- and $1$-forms satisfy the conditions (i)--(iii) in Theorem~\ref{T.main}, while the validity of the Theorem for $2$- and $3$-forms stems from the case of $0$- and $1$-forms and the commutativity of the Laplacian and the Hodge star operator.  
\end{proof}

Let us also notice the following corollary, which follows from the proof of Theorem~\ref{T.main} that we have presented above. 

\begin{corollary}\label{C.Beltrami}
For any $g\in\Ga$, the eigenvalues of the Beltrami operator $*_gd$ on $E_g$ are simple and the zeros of the corresponding eigenfunctions are all hyperbolic.
\end{corollary}

\section*{Acknowledgements}

This work has been supported in part by the Spanish Ministry of Science under grants no.\ FIS2008-00209 (A.E.) and MTM2007-62478 (D.P.-S.) and by Banco Santander--UCM under grant no.\ GR58/08-910556 (A.E.). The authors acknowledge the MICINN for financial support through a postdoctoral fellowship (A.E.) and the Ramón y Cajal program (D.P.-S.). The second author thanks the ETH Zürich for hospitality and support.

\end{document}